\documentclass[a4paper,oneside,english]{amsart}
\usepackage[T1]{fontenc}
\usepackage[latin9]{inputenc}
\usepackage{float}
\usepackage{textcomp}
\usepackage{graphicx}
\usepackage{setspace}
\usepackage{amssymb}

\makeatletter

\providecommand{\tabularnewline}{\\}

\numberwithin{equation}{section} 
\numberwithin{figure}{section} 
  \@ifundefined{theoremstyle}{\usepackage{amsthm}}{}
  \theoremstyle{plain}
  \newtheorem{thm}{Theorem}[section]
  \theoremstyle{definition}
  \newtheorem{defn}[thm]{Definition}
\newenvironment{lyxlist}[1]
{\begin{list}{}
{\settowidth{\labelwidth}{#1}
 \setlength{\leftmargin}{\labelwidth}
 \addtolength{\leftmargin}{\labelsep}
 }}
{\end{list}}
  \theoremstyle{remark}
  \newtheorem{rem}[thm]{Remark}
  \theoremstyle{plain}
  \newtheorem{algorithm}[thm]{Algorithm}

\makeatother

\usepackage{babel}

\begin{document}

\title{An algorithmic implementation of the $\pi$ function based on a new
sieve}

\author{Damián Gulich$^{A,B}$, Gustavo Funes$^{A,B}$, Nahuel Lofeudo$^{C}$,
Leopoldo Garavaglia$^{D}$, Mario Garavaglia$^{A,B}$}

\maketitle
\begin{singlespace}
\begin{center}
{\small $^{A}$Departamento de Física, Facultad de Ciencias Exactas,
Universidad Nacional de La Plata, Argentina }\\
{\small $^{B}$Laboratorio de Procesamiento Láser, Centro de Investigaciones
Ópticas, La Plata, Argentina }\\
{\small $^{C}$Facultad de Informática, Universidad Nacional de La
Plata, Argentina }\\
{\small $^{D}$Getafe, Madrid, Spain}
\par\end{center}{\small \par}
\end{singlespace}

\begin{singlespace}

\email{\emph{\small E-mail addresses:}{\small{} dgulich@ciop.unlp.edu.ar
(D. Gulich), gfunes@ciop.unlp.edu.ar (G. Funes), nlofeudo@lifia.info.unlp.edu.ar
(N. Lofeudo), garavaglia\_leo@hotmail.com (L.Garavaglia), garavagliam@ciop.unlp.edu.ar
(M. Garavaglia)}}
\end{singlespace}

\begin{abstract}
In this paper we propose an algorithm that correctly discards a set
of numbers (from a previously defined sieve) with an interval of integers.
Leopoldo's Theorem states that the remaining integer numbers will
generate and count the complete list of primes of absolute value greater
than 3 in the interval of interest. This algorithm avoids the problem
of generating large lists of numbers, and can be used to compute (even
in parallel) the prime counting function $\pi\left(h\right)$. 
\end{abstract}


\section{Introduction}

An old problem in mathematics is the way to compute the amount of
prime numbers less or equal to a given value $h$ \cite[p. 347]{Dickson}.
This function is known as $\pi\left(h\right)$ \cite{Weisstein}.
The preeminent method for such task since the 3rd century BC was the
sieve of Eratosthenes. From then on there were no great advances on
the subject until Gauss in 1863 \cite[p. 352]{Dickson}. His work
allowed more advances \cite{Hardy,Deleglise} and then boosted by
the growth of calculation power in the 20th century. Recent implementations
\cite{Deleglise} require several sophisticated computational strategies.

One of the major difficulties in solving this issue -or may be the
only one- is the question of how prime numbers are distributed among
integers. Matheaticians and other scientists have carried great efforts
to find some general formula and simple prime generating formula.
Such formula is still ellusive. Some formulas have been proposed,
but they are valid only in certain cases, such as Mersenne numbers,
and others. The key to find a simple formula that would allow us determine
all prime numbers, should begin demonstrating that all numbers obey
a general logic order, from there on it could be possible to find
a sequence that would permit us to establish such formula. In other
words it is about establishing a new sieve of universal validity.
This basic essential aspect was the one which guided L. Euler in 1770
to express \emph{\textquotedblleft{}Mathematicians have tried in vain
to discover some order in the sequence of prime numbers but we have
every reason to believe that there are some mysteries which the human
mind will never penetrate\textquotedblright{}}.

On the other side, the efforts of mathematicians like Mersenne, Fermat,
and many others from Ancient Greece up to our days, to introduce certain
logic in the calculus of some prime numbers with simple arithmetic
operations became useless due to the observation made by Gauss \cite{Dickson}
himself, who stated that the amount of prime numbers less than $h$
is near the value of the integral: \[
\mbox{Li}\left(h\right)=\int_{2}^{h}\frac{dt}{\ln\left(t\right)}\]

This observation allowed Gauss to introduce the function $\pi\left(h\right)$
toghether with a probabilistic argument in the studies about the distribution
of primes. In fact the function $\pi\left(h\right)$ is related with
the \emph{Fehlerintegral} or Error function (Erf), introduced by Gauss
to analyze the uncertainties that appear in experimental measurements.
It is to be noted, that in experimental sciences like Physics or Biology,
the scientist needs to justify that the results of a certain experiment
can be treated inside the gaussian frame. This a priori it does not
seem to happen in the mathematics area, in particular in Number Theory,
were Axioms play a fundamental role. Nevertheless, in the last years
a lot of works have merged regarding the issue of the location of
prime numbers, whether and if primes are located in a random way or
if they have a totally determined structure. The papers that were
published in this area are various and interesting (see references
\cite{wolf,wolf2,Bonnano,scafetta,Spriro2,Ares}). These results are
not completely satisfactory, therefore \textquotedblleft{}the human
mind\textquotedblright{} \textendash{}according to L. Euler\textendash{}
will continue the attempt to penetrate the mysteries of the apparent
unbreakable tower of primes. In particular, about this matter, Terence
Tao \cite{Terence tao} proposes the existence of certain dichotomy
between structure and randomness among primes, and tries to separate
 the structured component (or correlated) from the \textquotedblleft{}pseudorandom\textquotedblright{}
component (or decorrelated). Anyway, the existence of a distribution
of prime numbers affected by the dichotomy of structured-randomness
implies that the structured behavior or low complexity behavior can
be determined by any general formula, meanwhile in the high complexity
behavior we could recognize at least two aspects: i) A totally randomness
aspect or decorrelated, and the other, ii) A chaotic one, which obeys
some kind of logic formulation. In our opinion any well known prime
number allows us to describe a sieve, which can be enlarged as far
as new primes are verified as such numbers \cite{Elementary}. However
we consider that the unknown prime numbers, are not undefined nor
undetermined. That means that all of them can be defined and determined
by a general and unique algorithm, logically determined and located
sequentially among integers \cite{Garavaglia}.

Returning to our primary purpose of implementing an algorithm to calculate
the value of the $\pi(h)$ function for a given $h$, in this paper
we will study this problem defining a new sieve whose properties permit
an elementary study of $\pi(h)$, toghether with the possibility of
finding its value on a given interval. The basis of this sieve stands
on an infinite matrix presented in \cite{Elementary}, which has the
property of generating all indexes that lead to non prime numbers
according to the form $6n+1$ or $6n-1$, a similar approach   can
be found in \cite{Iovane}, but with a different implementation.

The paper is organized  as follow. In Section 2. we reviewed the properties
of the form $6n+1$, being $n=0,1,2,\ldots$ , which is the base of
our research. Section 3 is devoted to our proposal of a sieve to define
and compute all prime numbers larger than 3. In Sections 4 and 5 we
present the Leopoldo\textasciiacute{}s Theorem and the possibility
to compute  $\pi(h)$ function according with the $\Lambda$ algorithm
that we are introducing. Section 6 and 7 are devoted to our actual
results, comments on future research work and conclusions.

\section{About the form $6n+1$}

In \cite{Elementary} we reviewed some well known properties of numbers:

\begin{thm}
\label{thm:Modulos}Every prime number of absolute value greater than
3 can be written in the form $6n+1$ or $6n-1$.
\end{thm}
\begin{proof}
Let's see the equivalences modulo 6. Suppose $q$ prime.

1) If $q\underset{6}{=}0$ $\Rightarrow$$q=6n$$\Rightarrow$$6\mid q$,
ABS.

2) If $q\underset{6}{=}1$ $\Rightarrow$$q=6n+1$, which is not impossible
since $7=6+1$ is a prime.

3) If $q\underset{6}{=}2$ $\Rightarrow$$q=6n+2=2\left(3n+1\right)$$\Rightarrow$$2\mid q$,
ABS.

4) If $q\underset{6}{=}3$ $\Rightarrow$$q=6n+3=3\left(2n+1\right)$$\Rightarrow$$3\mid q$,
ABS.

5) If $q\underset{6}{=}4$ $\Rightarrow$$q=6n+4=2\left(3n+2\right)$$\Rightarrow$$2\mid q$,
ABS.

6) If $q\underset{6}{=}5$ $\Rightarrow$$q=6n+5=6n+6-1=6\left(n+1\right)-1$,
which is not impossible since with $n=1$ this gives 11, a prime.

Then we introduce some definitions:
\end{proof}
\begin{defn}
The $\alpha$ class of integer numbers \cite{Garavaglia} is the set
\begin{equation}
\alpha=\left\{ x\in\mathbb{Z}/x=6n+1,n\in\mathbb{Z}\right\} \label{eq:6n+1}\end{equation}

\end{defn}
~

\begin{defn}
The $\beta$ class of integer numbers \cite{Garavaglia} is the set

\begin{equation}
\beta=\left\{ x\in\mathbb{Z}/x=6n-1,n\in\mathbb{Z}\right\} \label{eq:6n-1}\end{equation}

\end{defn}
The different values of relations in \eqref{eq:6n+1} and \eqref{eq:6n-1}
are shown in Table \ref{tab:Valores-para-6n-1}.

\begin{center}
\begin{table}[H]
\begin{centering}
\begin{tabular}{|c|c|c|}
\hline 
$n$ & $\beta_{n}=6n-1$ & $\alpha_{n}=6n+1$\tabularnewline
\hline
\hline 
$\vdots$ & $\vdots$ & $\vdots$\tabularnewline
\hline 
-5 & -31 {*} & -29 {*}\tabularnewline
\hline 
-4 & -25 & -23{*}\tabularnewline
\hline 
-3 & -19 {*} & -17 {*}\tabularnewline
\hline 
-2 & -13 {*} & -11 {*}\tabularnewline
\hline 
-1 & -7 {*} & -5 {*}\tabularnewline
\hline 
0 & -1 & 1\tabularnewline
\hline 
1 & 5 {*} & 7 {*}\tabularnewline
\hline 
2 & 11 {*} & 13 {*}\tabularnewline
\hline 
3 & 17 {*} & 19 {*}\tabularnewline
\hline 
4 & 23 {*} & 25\tabularnewline
\hline 
5 & 29 {*} & 31 {*}\tabularnewline
\hline 
$\vdots$ & $\vdots$ & $\vdots$\tabularnewline
\hline
\end{tabular}
\par\end{centering}

\caption{\label{tab:Valores-para-6n-1}Values of $6n-1$ and $6n+1$ for several
$n$. With ({*}) we mark prime numbers.}

\end{table}

\par\end{center}

Strictly speaking, a {}``complete'' list of all prime numbers of
absolute value greater than 3 ($\left\{ \ldots,-7,-5,5,7,\ldots\right\} $)
is the list of primes from both classes, $\alpha$ and $\beta$.

We now state a property given in \cite{Garavaglia} and \cite{Elementary}: 

\begin{thm}
Every prime number of absolute value greater than 3 (except for the
sign) is generated by $6n+1$, with $n$ integer.
\end{thm}
\begin{proof}
We must prove the equivalence (except for the sign) between both families
given in Theorem \ref{thm:Modulos}. 

Let be $f_{\alpha}\left(n\right)=6n+1$ and $f_{\beta}\left(n\right)=6n-1$,
we must now prove that $f_{\alpha}\left(-n\right)=-f_{\beta}\left(n\right)$.
Indeed:\[
f_{\alpha}\left(-n\right)=6\left(-n\right)+1=-6n+1=-\left(6n-1\right)=-f_{\beta}\left(n\right)\]

\end{proof}
\begin{defn}
We define the set of integer numbers $G_{\alpha}$:\begin{equation}
G_{\alpha}=\left\{ g\in\mathbb{Z}/6g+1\mbox{ is a prime}\right\} \label{eq:G+}\end{equation}

This means, $G_{\alpha}$ is the set of \emph{all} numbers that (except
for the sign) generate \emph{all} primes of absolute value greater
than 3 by the relationship \eqref{eq:6n+1}.

In other words, this result is perfectly logic because it is based
on: 
\end{defn}
\begin{lyxlist}{00.00.0000}
\item [{i)}] Peano\textasciiacute{}s Axioma applied to both sequences of
numbers $\alpha$ and $\beta$ starting in 1 and 5, respectively,
taking into account that the next number in those sequences are defined
as 1 + 6, and 5 + 6, respectively, and so on and so forth; and 
\item [{ii)}] Remark \eqref{rem:Posibilidades} which follows from Theorem
3.4., represented by the Table of Products \ref{tab:Tabla-de-productos.}.
This table states that products between numbers belonging to the $\alpha$
sequence and products between numbers belonging to the $\beta$ sequence,
results in numbers of Class $\alpha$; while the crossed products
between numbers of the sequence $\alpha$ with numbers of the sequence
$\beta$, results in numbers of Class $\beta$. Besides, note that
crossed products between numbers belonging to the six classes we presented
in Theorem 2.1. never result in numbers belonging to classes $\alpha$
or $\beta$. Table \ref{tab:Numbers-belongning-to} shows the results
of these products from $n=0$ to $n=10$; remember, prime numbers
greater that 3 are only allocated in classes $\alpha$ or $\beta$.
\end{lyxlist}
\begin{defn}
\begin{table}[h]
\begin{centering}
\begin{tabular}{|c|c|c|c|c|c|c|}
\hline 
$n$ & $6n+1$($\alpha$) & $6n+2$ & $6n+3$ & $6n+4$ & $6n+5$ ($\beta$) & $6n+6$\tabularnewline
\hline
\hline 
0 & 1 & \textbf{2} & \textbf{3} & $4=2\times2$ & \textbf{5} & $6=2\times3$\tabularnewline
\hline 
1 & \textbf{7} & $8=2\times4$ & $9=3\times3$ & $10=2\times5$ & \textbf{11} & $12=2\times6$\tabularnewline
\hline 
2 & \textbf{13} & $14=2\times7$ & $15=3\times5$ & $16=2\times8$ & \textbf{17} & $18=2\times9$\tabularnewline
\hline 
3 & \textbf{19} & $20=2\times10$ & $21=3\times7$ & $22=2\times11$ & \textbf{23} & $24=2\times12$\tabularnewline
\hline 
4 & $25=5\times5$ & $26=2\times13$ & $27=3\times9$ & $28=2\times14$ & \textbf{29} & $30=2\times15$\tabularnewline
\hline 
5 & \textbf{31} & $32=2\times16$ & $33=3\times11$ & $34=2\times17$ & $35=5\times7$ & $36=2\times18$\tabularnewline
\hline 
6 & \textbf{37} & $38=2\times19$ & $39=3\times13$ & $40=2\times20$ & \textbf{41} & $42=2\times21$\tabularnewline
\hline 
7 & \textbf{43} & $44=2\times22$ & $45=3\times15$ & $46=2\times23$ & \textbf{47} & $48=2\times24$\tabularnewline
\hline 
8 & $49=7\times7$ & $50=2\times25$ & $51=3\times17$ & $52=2\times26$ & \textbf{53} & $54=2\times27$\tabularnewline
\hline 
9 & $55=5\times11$ & $56=2\times28$ & $57=3\times19$ & $58=2\times29$ & \textbf{59} & $60=2\times30$\tabularnewline
\hline 
10 & \textbf{61} & $62=2\times31$ & $63=3\times21$ & $64=2\times32$ & $65=5\times13$ & $66=2\times33$\tabularnewline
\hline 
... & ... & ... & ... & ... & ... & ...\tabularnewline
\hline
\end{tabular}
\par\end{centering}

\caption{\label{tab:Numbers-belongning-to}Numbers belongning to the six classes
we divided natural numbers from $n=0$ to $n=10$. Prime numbers are
shown in bold letters.}

\end{table}

\end{defn}

\section{A proposed    new sieve}

\begin{defn}
Let $A$ be an infinite matrix whose element $a\left(i,j\right)$
\footnote{Coordinates are in the Cartesian sense.%
} is \begin{equation}
a\left(i,j\right)=i+j\left(6i+1\right)\label{eq:a_ij}\end{equation}
where $i,j\,\in\mathbb{Z}$.

In Table \ref{tab:Elements-of-the} we show the elements in the central
part of the $A$ matrix. Note that axis numbers also match this representation.
\end{defn}
\begin{table}[H]
\[
\begin{array}{ccccccccccc}
 &  &  &  &  & \vdots\\
 & -96 & -71 & -46 & -21 & 4 & 29 & 54 & 79 & 104\\
 & -73 & -54 & -35 & -16 & 3 & 22 & 41 & 60\\
 & -50 & -37 & -24 & -11 & 2 & 15 & 28\\
 & -27 & -20 & -13 & -6 & 1 & 8\\
\cdots & -4 & -3 & -2 & -1 & 0 & 1 & 2 & 3 & 4 & \cdots\\
 & 19 & 14 & 9 & 4 & -1\\
 & 42 & 31 & 20 &  & -2\\
 & 65 & 48 &  &  & -3\\
 & 88 &  &  &  & -4\\
 &  &  &  &  & \vdots\end{array}\]

\caption{\label{tab:Elements-of-the}Elements of the central portion of $A$}

\end{table}

\subsection{Properties}

\begin{thm}
\label{thm:A-es-sim=0000E9trica.}$A$ is symmetrical.
\end{thm}
\begin{proof}
A simple expansion shows that \[
a\left(i,j\right)=i+j\left(6i+1\right)=i+6ij+j=j+i\left(6j+1\right)=a\left(j,i\right)\]

\end{proof}
\begin{defn}
Let $\widetilde{A}$ be the set of unrepeated elements of $A$ excluding
the axes (the elements of the form $a\left(i,0\right)$ y $a\left(0,j\right)$).
\end{defn}

\subsubsection{About the signs of the elements of $\tilde{A}$}

\begin{figure}[H]
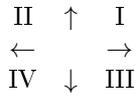

\[
\begin{array}{ccc}
\mbox{II} & \uparrow & \mbox{I}\\
\leftarrow &  & \rightarrow\\
\mbox{IV} & \downarrow & \mbox{III}\end{array}\]

\caption{\label{fig:Quadrants-of-A}Quadrants of $\tilde{A}$ }

\end{figure}

Four quadrants can be distinguished as shown in Figure \ref{fig:Quadrants-of-A}

What happens to the signs of the elements of $\tilde{A}$ from each
quadrant?

Because of Theorem \ref{thm:A-es-sim=0000E9trica.}, we should only
focus on the sign of elements of $\tilde{A}$ originally from quadrants
I, II, and IV.

\begin{itemize}
\item In quadrant I ($i\geq1$, $j\geq1$) all elements are positive
\item In quadrant II ($i\leq-1$, $j\geq1$)\\
It's easy to see that $j\left(6i+1\right)\leq0$, then $\tilde{a}\left(i,j\right)\leq0$
$\forall i,j$. 
\item In quadrant IV ($i\leq-1$, $j\leq-1$)\\
$i+j\left(6i+1\right)=i+j+6ij$. $\left(i+j\right)\leq-1$ y $ij\geq\left|i+j\right|\geq1$.
Then, the sign is positive.
\end{itemize}
\begin{thm}
\label{thm:Los-elementos-de-A-tilde-NO}The elements of $\widetilde{A}$
DO NOT generate prime numbers.
\end{thm}
\begin{proof}
$\tilde{a}\left(i,j\right)=i+j\left(6i+1\right)$ with $i\neq0$ and
$j\neq0$. If we put this into \eqref{eq:6n+1} and suppose $p$ prime,
then\[
p=6\tilde{a}\left(i,j\right)+1=6\left(i+j\left(6i+1\right)\right)+1=6i+6j\left(6i+1\right)+1=\left(6i+1\right)+6j\left(6i+1\right)=\left(6i+1\right)\left(6j+1\right)\]
 but since $i$ and $j$ are different from zero, then $p$ would
be a composite, ABS.
\end{proof}
\begin{rem}
\label{rem:Posibilidades}According to the signs of $i$ and $j$,
$6\tilde{a}\left(i,j\right)+1$ sweeps (except for the sign) all possibilities:
\end{rem}
\begin{enumerate}
\item If $i>0$ y $j>0$, the generated number is of the form $\alpha\cdot\alpha$.
\item If $i<0$ y $j<0$, the generated number is of the form $\beta\cdot\beta$.
\item If $i$ y $j$ have opposite signs, the generated number is of the
form $\alpha\cdot\beta$.
\end{enumerate}
See Table \ref{tab:Tabla-de-productos.} for properties of products
of $\alpha$'s and $\beta$'s.

\begin{center}
\begin{table}[H]
\begin{tabular}{|c||c|c|}
\hline 
$\times$ & $\alpha$ & $\beta$\tabularnewline
\hline
\hline 
$\alpha$ & $\alpha$ & $\beta$\tabularnewline
\hline 
$\beta$ & $\beta$ & $\alpha$\tabularnewline
\hline
\end{tabular}

\caption{\label{tab:Tabla-de-productos.}Table of products.}

\end{table}

\par\end{center}

Then, summarizing previous results, we can state that for positive
prime numbers:

\begin{enumerate}
\item The only even prime number 2 belongs to the sequence defined by the
form $2+6n$ for $n=0$. All the other even numbers in the same sequence
are multiples of 2, then they are composites.
\item All the others even numbers that belong to sequences of the form $4+6n$
and $6+6n$, being $n=0,1,2,...$, are multiples of 2 and obviously
they are composites.
\item The first odd prime number 3 belongs to the sequence defined by the
form $3+6n$ for $n=0$. All the other odd numbers in the same sequence
are multiples of 3, then they are composites.
\item Then, except 2 and 3, all the other prime numbers belong to classes
$\alpha$ ($1+6n$) or $\beta$ ($5+6n$).
\item Therefore, taking into account Table of Products $\alpha\times\alpha=\alpha$,
$\beta\times\beta$, and $\alpha\times\beta=\beta$, all composite
numbers generated by these products belong to classes $\alpha$ and
$\beta$.
\item Then all $\alpha$ prime numbers are those expressed by ($1+6n$)
for $n=0,1,2,...$, except:

\begin{enumerate}
\item Number 1 expressed by $1+6n$, being $n=0$.
\item All even and odd powers of all class $\alpha$ numbers, as $(1+6n)^{2+k}$
, for $n=1,2,3,...$, and $k=0,1,2,...$ 
\item All products between $\beta$ numbers and all their even and odd powers,
as: $[(5+6n)(5+6(n+m))]^{1+k}$, for $n=0,1,2,...$, $m=0,1,2,...$,
and $k=0,1,2,...$ 
\end{enumerate}
\item All $\beta$ prime numbers, except all the products of each $\beta$
number with all $\alpha$ numbers, except 1, like $(5+6n)(1+6m)$,
for $n=0,1,2,...$, and $m=1,2,3,4,...$ 
\end{enumerate}

\section{Leopoldo's Theorem and the $\pi$ function}

\subsection{Leopoldo's Theorem}

We defined the $\tilde{A}$ set as a list of all the non repeated
off-axis elements of $A$. A simple expansion sowed that the elements
of $\tilde{A}$ do not generate prime numbers by \eqref{eq:6n+1}.
Finally, we stated and proved Leopoldo's theorem \cite{Elementary}:

\begin{thm}
\label{thm:Teorema-de-Leopoldo}(Leopoldo's Theorem%
\footnote{We named theorem 4.1 as Leopoldo`s theorem because it was stated by
Leopoldo Garavaglia in the Summer 2006, Spain, but the basis was not
published until 2007 \cite{Garavaglia} and it was proved in \cite{Elementary}%
}) $G_{\alpha}=\mathbb{Z}-\widetilde{A}$.
\end{thm}
\begin{proof}
Every prime number of absolute value greater than 3 is either $\alpha$
or $\beta$, and products between those classes of equivalence are
closed on themselves. Because of Theorem \ref{thm:Los-elementos-de-A-tilde-NO}
and Remark \ref{rem:Posibilidades}, we know that $\widetilde{A}$
generates all possible $\alpha$ and $\beta$ composite numbers, and
so the elements $n$ of $\mathbb{Z}\notin\widetilde{A}$ generate
prime numbers (except for the sing) by $6n+1$.
\end{proof}
This means that all integers not generated by \eqref{eq:a_ij} will
generate all primes of absolute value greater than 3 by \eqref{eq:6n+1}
(and thus, this works as a sieve). 

In order to see the whole picture we are going to suppose that instead
of integers we are using real numbers for $i$ and $j$. This will
lead us to a graphical sight of the matrix where, keeping $a\left(i,j\right)$
as an integer, we should see hyperbolas. The plots will represent
any integer $c$ generated by the matrix, and they can be considered
as level sets. These sets have certain coincidences when $i$ and
$j$ are integers at the same time, obviously this is a point in which
the elements in the matrix coincide (in space) with an element in
the set of curves. Some of the curves have lots of repetitions specially
when $c$ is large, that is why the calculation of all the elements
of $A$ is futile. What we are looking here, in order to find primes,
are the integers $c$ which have integer generators $i$ and $j$.
This will lead us to composite numbers and indirectly to primes. Some
level sets of $A$ are shown in Figure \ref{eq:Curva-de-nivel}.

Now, suppose you wish to calculate all prime numbers of absolute value
greater than 3 up to a certain value $h=6c+1$ ($c>0$), this means
computing $\pi\left(h\right)$ using Leopoldo's Theorem. At first
sight, one would have to: 

\begin{enumerate}
\item Generate all elements of $\widetilde{A}$ up to $c$ which means that
$\left|a\left(i,j\right)\right|=\left(h-1\right)/6$
\item Discard the axis elements
\item Sort the rest
\item Discard repetitions
\item Remove them from the interval $\left[-c,c\right]$
\item Count the remaining numbers
\item Apply \eqref{eq:6n+1} to show the primes in the interval
\end{enumerate}
With large numbers, this computation would quickly become time and
memory prohibitive. However, a closer look at the distribution of
elements in the sieve, gives an appropriate answer to this problem.

\section{Leopoldo's Theorem and the computation of the $\pi$ function}

This method is based on discarding non prime generators. In fact we
only need to find two integers $i$ and $j$ to verify that the number
$c$ that we are checking belongs to $\widetilde{A}$. Once we find
these numbers there is no need for further search. This will simplify
the problem by avoiding repetitions. In this section we are going
to present an approach for an algorithmic method to compute $\pi$
based in the ideas presented above.

\begin{center}
\begin{figure}[H]
\begin{centering}
\includegraphics[width=8cm]{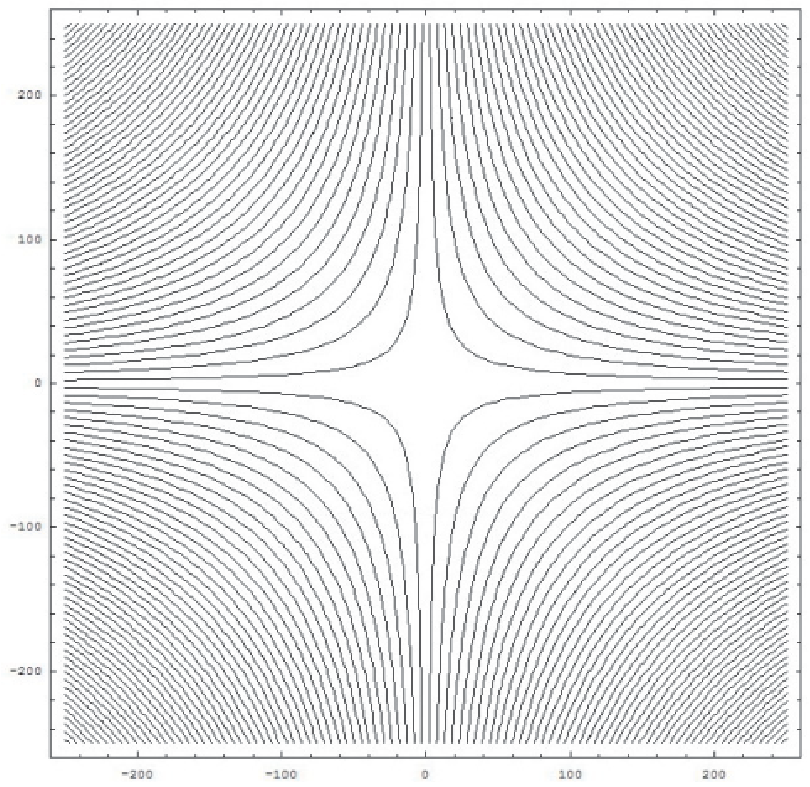}
\par\end{centering}

\caption{\label{fig:Curvas-de-nivel}Several level sets of $f\left(x,y\right)=x+y\left(6x+1\right)$. }

\end{figure}

\par\end{center}

\begin{center}
\begin{figure}[H]
\begin{centering}
\includegraphics[width=12cm]{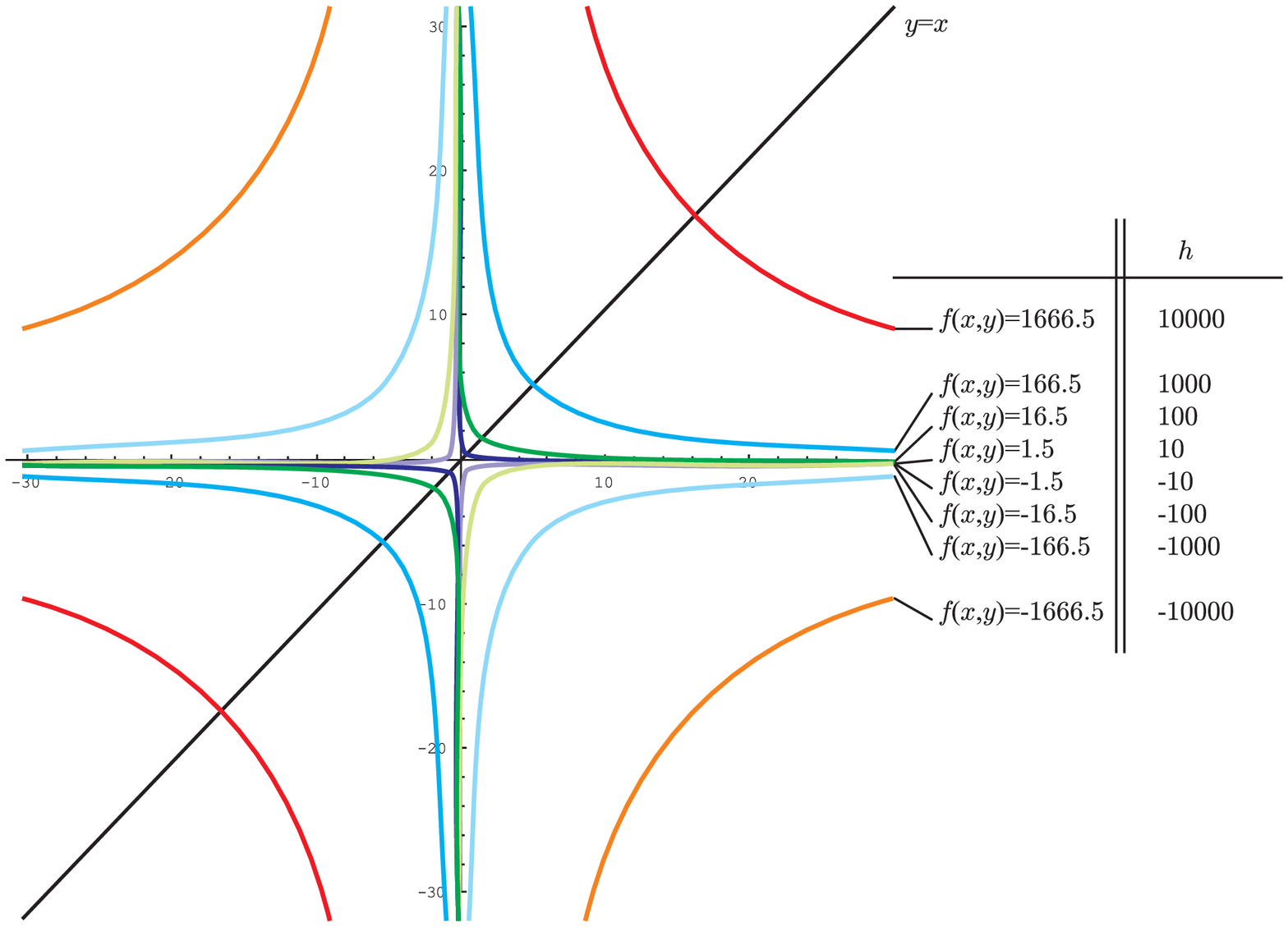}
\par\end{centering}

\caption{\label{fig:Representaci=0000F3n-de-las}Level sets of $f\left(x,y\right)=x+y\left(6x+1\right)=\pm c$. }

\end{figure}

\par\end{center}

In Figure \ref{fig:Representaci=0000F3n-de-las}, we show the representation
of the level sets $f\left(x,y\right)=\pm c$ ($c>0$). We must only
consider non-repeated elements of $\tilde{A}$ originally from within
the {}``star'' delimited by \begin{equation}
f\left(x,y\right)=\pm c=\pm\left(h-1\right)/6\label{eq:Curva-de-nivel}\end{equation}

\subsection{An algorithmic approach}

An exploration of the level sets allows algorithmic approach to $\pi\left(h\right)$.

\begin{algorithm}
We define the $\Lambda$ algorithm of arguments $c_{1}$ and $c_{2}$
($c_{2}>c_{1}\geq8$), as the procedure that
\end{algorithm}
\begin{enumerate}
\item Declares a natural variable $L=0$
\item For $c$ taking every integer value from $c_{1}$ to $c_{2}$

\begin{enumerate}
\item \label{enu:paso_1}For integers $x$ from $x=-\left\lfloor \left(c+1\right)/5\right\rfloor $
to $x=-\left\lfloor \left(\sqrt{1+6c}+1\right)/6\right\rfloor $ verifies
if \[
\frac{c-x}{6x+1}\]
 takes an integer value%
\footnote{This means, $\left(c-x\right)\equiv0$ (modulo $6x+1$).%
}. 

\begin{enumerate}
\item if it finds one, adds 1 to $L$ and goes to step (2)(c)
\item if it doesn't find any, goes to step (2)(b)
\end{enumerate}
\item \label{enu:Paso_2}For integers $x$ from $x=1$ to $x=\left\lfloor \left(\sqrt{1+6c}-1\right)/6\right\rfloor $
verifies if \[
\frac{c-x}{6x+1}\]
 takes an integer value.

\begin{enumerate}
\item if it finds one, adds 1 to $L$ and goes to step (2)(c)
\item if it doesn't, prints $c$ and $6c+1$, and goes to step (2)(c)
\end{enumerate}
\item \label{enu:Paso_3}For integers from $x=-\left\lfloor \left(c+1\right)/7\right\rfloor $
to $x=-1$ verifies if \[
\frac{-c-x}{6x+1}\]
 takes integer values.

\begin{enumerate}
\item if it finds one, adds 1 to $L$ and goes to the next value of $c$.
\item if it doesn't find any, prints $-c$ and $-6c+1$, and goes to the
next value of $c$.
\end{enumerate}
\end{enumerate}
\item Once the process has been completed up to $c_{2}$, reports the accumulated
value of $L$: \[
\Lambda\left(c_{1},c_{2}\right)=L_{final}\]

\end{enumerate}
This algorithm tells the amount of numbers that \emph{do not} generate
prime numbers by $6c+1$ in the interval $\left[c_{1},c_{2}\right]$,
and also gives the list of the missing values as well as the primes
generated by them. So, the amount of prime numbers between $h_{1}=6c_{1}-1$
and $h_{2}=6c_{2}+1$ ($c_{1}>8$) is: \begin{equation}
\Delta\pi=2\left(c_{2}-c_{1}\right)-\Lambda\left(c_{1},c_{2}\right)+1\label{eq:P}\end{equation}
In the particular case of $c_{1}=8$ ($h_{1}=47$) and $h=6c_{2}+1$:\begin{equation}
\pi\left(h\right)=2c_{2}-\Lambda\left(8,c_{2}\right)\label{eq:pi_exacto}\end{equation}

\section{Results and Future Works}

In Table \ref{tab:Resultados}, it may be seen that the the proposed
algorithm for the $\pi$ function agrees with known results for several
testing values. Equations \eqref{eq:P} and \eqref{eq:pi_exacto}
enable a parallelization of the computation of the $\pi$ function
with the $\Lambda$ algorithm.

Calculation time grows with the parameter $c$. The last value in
the table took nearly an hour to be computed with an AMD Athlon 64
X2 Dual Core 4200+. All values were calculated using \eqref{eq:pi_exacto}.
The memory used remains stable, since it doesn't take more than it
needs to store the variables and operations in steps (2)(a), (2)(b)
and (2)(c). 

\begin{center}
\begin{table}[!h]
\begin{tabular}{|c||c|c|c|}
\hline 
$h$ & $\pi\left(h\right)$ & $\pi\left(h\right)$ & $\pi\left(h-3\right)$\tabularnewline
 & (Mathematica 5) & ($\Lambda$ Algorithm) & \cite{Weisstein}\tabularnewline
\hline
\hline 
$10^{2}+3$ & 27 & 27 & 25\tabularnewline
\hline 
$10^{3}+3$ & 168 & 168 & 168\tabularnewline
\hline 
$10^{4}+3$ & 1229 & 1229 & 1229\tabularnewline
\hline 
$10^{5}+3$ & 9593 & 9593 & 9592\tabularnewline
\hline 
$10^{6}+3$ & 78499 & 78499 & 78498\tabularnewline
\hline 
$10^{7}+3$ & 664579 & 664579 & 664579\tabularnewline
\hline
\end{tabular}

\caption{\label{tab:Resultados}$\Lambda$ algorithm based results versus several
known values of $\pi\left(h\right)$.}

\end{table}

\par\end{center}

Given the fact that this algorithm evaluates all values between $c_{1}$
and $c_{2}$, as well as the independence of the obtained value $\Delta\pi$
with other intervals, a parallel expression of the Algorithm is possible
using $n$ independent processors where the range of values $c$ for
each acting processor is\[
C_{1P}=c_{1}+\left(\frac{c_{2}-c_{1}}{\#P}\right)P_{n}\]
\[
C_{2P}=C_{1P}+\left(\frac{c_{2}-c_{1}}{\#P}\right)P_{n}\]
 and so on, where $\#P$ is the total number of processors and $P_{n}$
is the number of a given processor in the cluster.

Since the computing load is equally divided between computing nodes,
the efficiency of the process asymptotically approaches 1 (each computing
node uses nearly 100\% of itself to calculate). The speedup tends
to $\#P$ as new computing nodes are added whenever $\left(c_{2}-c_{1}\right)\gg\#P$.
In a future paper we will deal with the case $\left(c_{2}-c_{1}\right)\leq\#P$

\section{Conclusions}

We have seen that using a simple classification of numbers we can
group primes in families. We have also seen that primes and its multiples
are grouped in the same classes according to the table of products
\ref{tab:Tabla-de-productos.}. This classification leads us to the
separation of composite numbers among primes, and the characteristic
index number that do not generate primes. The matrix $A$ which is
the matrix of indexes (without the axis) and we use it in order to
find primes. Also by analyzing  the matrix as continuous hyperbolas
we have found a way to count the amount of elements of $A$ from a
set of integers. In this way by counting the amount of {}``non prime
generating indexes'' we can obtain indirectly the amount of primes
in a given section. The advantage of this method is that we can calculate
the amount of primes by section. That means that this is not only
an accumulative method, we can use it in any section of natural numbers
and also use various computers to process different sections at the
same time. Finally we can also use known results of $\pi(h)$ and
start counting primes from the last well known result. As mentioned
above, in future works we would like to analyze the use of several
computers in order to find $\pi(h)$ and test the proposed algorithm.

\section{Acknowledgments}

Damián Gulich and Gustavo Funes are financially supported by a student
fellowship from the INNOVATEC Foundation, Argentina.

Damián Gulich and Gustavo Funes thank Dr. Mario Garavaglia for involving
them in this line of research.


\begin{thebibliography}{10}
\bibitem{Elementary}Damian Gulich, Gustavo Funes, Leopoldo Garavaglia,
Beatriz Ruiz, Mario Garavaglia, (2007), ''An elementary sieve''.
arXiv:0708.3709v1 {[}math.GM]. \texttt{http://arxiv.org}

\bibitem{Dickson}Dickson, Leonard Eugene. (1952), \emph{History of
the theory of numbers}, (Vol. 1), New York, N. Y.: Chelsea Publishing
Company.

\bibitem{Hardy}Hardy, G. H y Wright, E. M. (1962), \emph{An introduction
to the theory of numbers}, (4th ed.), Oxford: Oxford at the Clarendon
Press.

\bibitem{Weisstein}Weisstein, Eric W. \char`\"{}Prime Counting Function.\char`\"{}
From \emph{MathWorld}--A Wolfram Web Resource. \texttt{http://mathworld.wolfram.com/PrimeCountingFunction.html}

\bibitem{Garavaglia}Garavaglia, Leopoldo and Garavaglia, Mario. (2007),
{}``On the location and classification of all prime numbers''. arXiv:0707.1041v1
{[}math.GM]. \texttt{http://arxiv.org}

\bibitem{Deleglise}Deleglise, M. and Rivat, J. (1996), {}``COMPUTING
$\pi(x)$: The Meissel, Lehmer, Lagarias, Miller, Odlyzko method''.
Mathematics of computation (Vol. 65, Nº 213). Jan. 1996, P. 235-245.

\bibitem[7]{Iovane}Gerardo Iovane (Salerno U. \& Frascati) (2008),
{}``The set of prime numbers: Symmetries and supersymmetries of selection
rules and asymptotic behaviours''. 12pp. Chaos Solitons Fractals
37: 950-961,2008.

\bibitem[8]{Terence tao}Terence Tao, (2005), {}``The dichotomy between
strucure and randomness, arithmetic progressions, and the primes''.
arXiv:math/0512114v2 {[}math.NT] \texttt{http://arxiv.org}

\bibitem[9]{wolf} Mark Wolf, (1997), {}``1/f noise in the distribution
of prime numbers''. Physics A, 241: 493-499, 1997.

\bibitem[10]{wolf2} Mark Wolf, (1998), {}``Random walk on the prime
numbers''. Physics A, 250: 335-344, 1998.

\bibitem[11]{Bonnano} Bonanno Claudio, Mirko S Mega, (2004), {}``Toward
a dinamical model of prime numbers''. Chaos Solitons \& Fractals
20: 107-118, 2004

\bibitem[12]{Stanley} Pradep Kumar, Plamen Ch Ivanov, Stanley H Eugene,
(2008), {}``Information entropy and correlations in prime numbers''.
arXiv:cond-mat/0303110v4 {[}cond-mat.stat-mech] \texttt{http://arxiv.org}

\bibitem[13]{scafetta} N Scafetta, T Imholt, J A Roberts, B J West,
(2004), {}``An intensity-expansion method to treat non-stationary
time series: an application to the distance between prime numbers''.
Chaos Solitons \& Fractals 20: 119-125, 2004.

\bibitem[14]{Ares} S Ares, M Castro, (2006), {}``Hidden structure
in a randomness of the prime number sequence?''. Physica A 360: 285-296,
2006.

\bibitem[15]{Spiro} George G Szpiro, (2007), {}``Peacks and gaps:
spectral analysis of the intervals between prime numbers''. Physica
A, 384: 291-296, 2007.

\bibitem[16]{Spriro2}George G Szpiro, (2004), {}``The gaps between
the gaps: some patterns in the prime number sequence''. Physica A,
341: 607-617, 2004.
\end{thebibliography}
\end{document}